\documentclass[symmetry,article,accept,moreauthors,pdftex,a4paper]{mdpi} 
\usepackage[english]{babel}
\usepackage{amsmath,xparse}  
\usepackage{amsfonts}
\usepackage{amssymb}  
\usepackage{mathrsfs} 
\usepackage{graphicx}
\usepackage{color}
\usepackage{caption}
\usepackage{subcaption}
\usepackage{mathtools}
\usepackage{graphicx}
\usepackage{kotex}

\firstpage{1} 
\pubvolume{14}
\issuenum{1}
\articlenumber{54}
\pubyear{2022}
\copyrightyear{2022}
\externaleditor{Academic Editor: Danny Arrigo} 
\datereceived{26 November 2021} 
\dateaccepted{19 December 2021} 
\datepublished{1 January 2022} 
\hreflink{https://doi.org/\linebreak10.3390/sym14010054} 

\Title{On Modified Second Paine--de Hoog--Anderssen Boundary Value Problem}

\TitleCitation{On Modified Second Paine--de Hoog--Anderssen Boundary Value Problem}

\Author{Natanael Karjanto \orcidA{}}

\AuthorNames{Natanael Karjanto}

\AuthorCitation{Karjanto, N.}

\address[1]{%
Department of Mathematics, University College, Natural Science Campus,  Sungkyunkwan University, Suwon~16419, Korea; natanael@skku.edu\\
}

\abstract{This article deals with a special case of the Sturm--Liouville boundary value problem (BVP), an eigenvalue problem characterized by the Sturm--Liouville differential operator with unknown spectra and the associated eigenfunctions. By examining the BVP in the Schr\"odinger form, we are interested in the problem where the corresponding invariant function takes the form of a reciprocal quadratic form. We call this BVP the modified second Paine--de Hoog--Anderssen (PdHA) problem. We estimate the lowest-order eigenvalue without solving the eigenvalue problem but by utilizing the localized landscape and effective potential functions instead. While for particular combinations of parameter values that the spectrum estimates exhibit a poor quality, the outcomes are generally acceptable although they overestimate the numerical computations. Qualitatively, the eigenvalue estimate is strikingly excellent, and the proposal can be adopted to other BVPs.}


\begin{document}

\section{Introduction}

The Sturm--Liouville boundary value problem (BVP) is an active research area in mathematics and its applications~\cite{kravchenko2020direct,zettl2021recent}. While the formulation originally comes from applied problems~\cite{gwaiz2008sturm,haberman2013applied}, the mathematical theory has reached maturity, both in terms of analysis {and} its numerical solutions~\cite{teschl2021ordinary,zettl2010sturm,amerin2005sturm,pryce1993numerical,bailey2001sleighn2}. 

One of the intriguing problems in the Sturm--Liouville BVP is finding eigenvalues and their corresponding eigenfunctions. Some attempts to estimate eigenvalues were performed more than a half-century ago~\cite{hochstadt1961,pruess1973estimating}. Various numerical techniques have been developed, including the Rayleigh--Ritz method~\cite{trefethen1997numerical}, matrix-variational method~\cite{gerck1982solution}, shooting method~\cite{ledoux2009efficient}, Pruess' coefficient approximation method~\cite{pruess1975high}, spectral parameter power series method~\cite{kravchenko2010spectral} and modified integral series methods among others~\cite{moan1988efficient}. A body of literature is continuously being updated by the modification and improvement of the existing techniques, as well as with applications of the Sturm--Liouville problem in numerous physical situations.

The history of the Sturm--Liouville problem goes back to the 19th century when the mathematicians Jacques Charles Fran\c{c}ois Sturm (1803--1855) from Geneva, then part of France, and Joseph Liouville (1809--1882) from France, investigated particular problems of second-order linear differential equations under appropriate boundary conditions and the properties of the corresponding solutions. Their results were published in a sequence of papers in 1836--1837. 

The problems arise frequently not only as governing equations of motion based on Newtonian mechanics but also when dealing with the method of separation of variables in tackling linear partial differential equations (PDEs). There are also applications in quantum mechanics where the time-independent Schr\"odinger equation for a particle with fixed angular momentum quantum numbers moves in a spherically symmetric potential at the energy level associated with spectra~\cite{prugovecki1981quantum,newton2013scattering}.

The classical Sturm--Liouville problem is composed of a linear second-order ordinary differential equation (ODE) and a set of boundary conditions at the endpoints of the interval where the problem is defined. It can be written in the following \emph{self-adjoint} (or \mbox{\emph{canonical}) form}
\begin{align}
-\frac{d}{dz} \left(u(z) \frac{dy}{dz}\right) + v(z) y &= \lambda w(z) y, \qquad \qquad z_0 < z < z_1, 	\label{ODE}\\
a_0 y(z_0) - a_1 u(z_0) \frac{dy}{dz}(z_0) = 0,& \qquad \qquad b_0 y(z_1) + b_1 u(z_0) \frac{dy}{dz}(z_1) = 0.	\label{BC}
\end{align}

An \emph{eigenvalue} (or \emph{spectrum}) is a value of $\lambda$ such that the ODE~\eqref{ODE} possesses a nontrivial solution $y(z)$ subject to the prescribed boundary conditions. This solution is called the corresponding \emph{eigenfunction} associated with each $\lambda$ and is unique up to scalar multiplications. Note that while some authors impose the classical derivatives~$y'$ at the boundaries, we opt for quasi-derivatives~$(uy')$ instead, where prime denotes the differentiation with respect to the variable $z$. Furthermore, neither $a_0$ and $a_1$ nor $b_0$ and $b_1$ are both zero.

The \emph{regular} Sturm--Liouville problem operates under the assumption that the coefficient functions are well-behaved. In this case, $u(z)$ and $w(z)$ are strictly positive, and $u(z)$, $u'(z)$, $v(z)$, and $w(z)$ are continuous functions on the bounded closed interval $[z_0, z_1]$. As a consequence, there exists an infinite sequence of eigenvalues $\lambda_0 < \lambda_1 < \lambda_2 < \cdots$ and eigenfunctions $y_0(z)$, $y_1(z)$, $y_2(z)$, $\cdots$ such that $y_n(z)$ has only $n$ zeros on the open interval $(z_0, z_1)$. In addition, distinct eigenfunctions are orthogonal with respect to the weight function $w(z)$, i.e.,
\begin{equation*}
\int_{z_0}^{z_1} w(z) \, y_i(z) \, y_j(z) \, dz = 0, \qquad \qquad \text{whenever} \qquad i \neq j.
\end{equation*}

The \emph{singular} Sturm--Liouville problem comprises whether the coefficient functions have singularities at the boundary points or whether the interval is unbounded. While the regular Sturm--Liouville problem was introduced by and the focus of Sturm's and Liouville's works, the singular Sturm--Liouville problem was initiated by Hermann Weyl from G\"ottingen, Germany, who investigated some ODEs with singularities and proposed the topic of the essential spectrum in 1910~\cite{weyl1910ueber}. Indeed, the progressive development of quantum theory in the 1920s and 1930s and a breakthrough in the general spectral theorem for unbounded self-adjoint operators in the Hilbert space provided an impetus for further examination into the spectral theory of Sturm--Liouville self-adjoint differential operators~\cite{zettl2010sturm,titchmarsh1962eigenfunction}.

The body of literature on Sturm--Liouville theory is overwhelmingly plentiful. Thanks to the advancement in computational software and hardware, the numerical aspect of the Sturm--Liouville problem is particularly one of the most active research areas of mathematical physics. The following brief literature review is by no means exhaustive. While a comparison between the Sinc--Galerkin technique and differential transform method reveals that the latter is more efficient than the former~\cite{alquran2010approximations}, a comparative study between the Sinc--Galerkin technique and variational iteration method indicates that the former is better than the latter when dealing with singular problems~\cite{alkhaled2020comparison}.

During the past decade, both regular and singular fractional derivatives Sturm--Liouville problems and operators have gained a lot of attention~\cite{almdallal2009efficient,klimet2013fractional,zayernouri2013fractional}. In particular, Klimet and Agrawal showed for the first time that the eigenvalues of two classes of the fractional Sturm--Liouville operators are real-valued, and the eigenfunctions associated with two distinct eigenvalues are orthogonal~\cite{klimet2013fractional}. Furthermore, the existence of a countable set of orthogonal eigenfunctions for the regular fractional Sturm--Liouville problem was proven by the methods of fractional variational analysis~\cite{klimek2014variational}. Recently, an efficient numerical method for estimating eigenvalues and eigenfunctions of the Caputo-type fractional Sturm--Liouville problem based on the Lagrange polynomial interpolation was proposed~\cite{sadabad2021efficient}.
	
When it comes to applications, the Sturm--Liouville problem features an abundance of them in applied mathematics and physics. The BVP itself originally sprang from the heat conduction problem in a nonhomogeneous thin bar, but equally classic problems from PDE include the vibration of plucked strings, thin membranes, and sound waves. The Sturm--Liouville BVP emerges naturally as an immediate consequence of implementing the method of separation of variables. For these vibration problems, the eigenvalues determine the frequency of oscillation, whereas the associated eigenfunctions correspond to the shape of the vibrating waves at any point in time~\cite{luetzen1990joseph}.

The most straightforward applications of the Sturm--Liouville problem bring about the various Fourier series, and more sophisticated applications lead to the generalizations of Fourier series involving Bessel functions, Hermite polynomials, and other special functions. One example is an electrostatic field generated by a spherical capacitor, where Laplace's equation in spherical coordinates serves as the governing equation and the Legendre polynomials appear in the solution~\cite{gwaiz2008sturm}. Another example comes from a quantum-mechanical system, where the wave function of a particle is governed by a one-dimensional time-independent Schr\"odinger equation. In this case, the eigenvalues represent the energy levels of the atomic system and the eigenfunctions are the wave functions of observable physical quantities~\cite{band2013quantum,griffiths2018introduction}.

The focal point of this article is the so-called \emph{modified second Paine--de Hoog--Anderssen (PdHA) BVP}~\cite{paine1981correction}. What we refer to as the ``first'' and ``second'' PdHA problems deal with the corresponding invariant function of the ODE in the canonical form in the Sturm--Liouville---in this case, the PdHA problem. While the first PdHA problem employs an exponential invariant function, i.e., $q(\widehat{z}) = e^{\widehat{z}}$, the second ``classical'' PdHA problem adopts a special case of the reciprocal quadratic binomial function, i.e., $q(\widetilde{z}) = (\widehat{z} + 0.1)^{-2}$. We do not examine the former in this article. Our focus is to investigate the latter by modifying the invariant function to involve several parameters, i.e., $q(\widetilde{z}) = c \, (a\widehat{z} + b)^{-2}$, where $a, b$ and $c > 0$. Hence, the name of the \emph{modified second PdHA problem}. Even though the classical problem appeared more than four decades ago, investigating the modified version still proves to be illuminating, as we will encounter in this article.

Since finding eigenvalues and eigenfunctions is a central aspect of the Sturm--Liouville problem, our motivation is to provide an alternative for estimating the former without actually solving the original modified second PdHA BVP. Our main contribution is a quantitative comparison between the eigenvalues obtained from an estimate with numerical solution results. By combining different parameters, we observe how the eigenvalues behave as we vary the parameters. Although the particular choice of the problem is rather narrow and specific, we anticipate that our contribution may illuminate other problems that can be solved using a similar technique. We also invite comments and debates from expertise for potential approaches to improve the estimated accuracy.

This article is organized as follows. Section~\ref{general} deals with the modified second PdHA problem. Using a transformation well-known in the Sturm--Liouville theory, we verify that a particular Sturm--Liouville problem in the canonical form can be transformed into the modified second PdHA BVP, where the corresponding invariant function admits the form of reciprocal quadratic binomial function. Section~\ref{estimate} covers how to estimate the spectrum without solving an eigenvalue problem but by incorporating the landscape function and effective potential instead. Section~\ref{compare} compares the eigenvalue estimate with numerical simulations for both Dirichlet and Neumann boundary conditions. Finally, Section~\ref{conclude} concludes our discussion.  
 
\section{Modified Second Paine--de Hoog--Anderssen Problem}	\label{general}

\begin{figure}[H]
\widefigure
\begin{subfigure}{0.45\textwidth}
\caption{} \vspace{0.1cm}
\includegraphics[width=\linewidth]{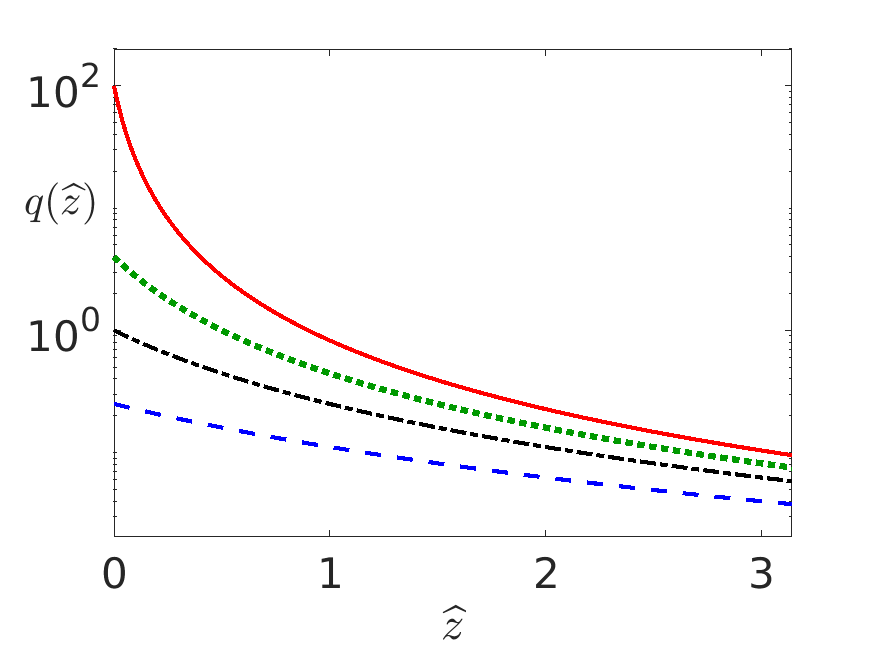}
\end{subfigure} \hspace*{0.5cm}
\begin{subfigure}{0.45\textwidth}
\caption{} \vspace{0.1cm}
\includegraphics[width=\linewidth]{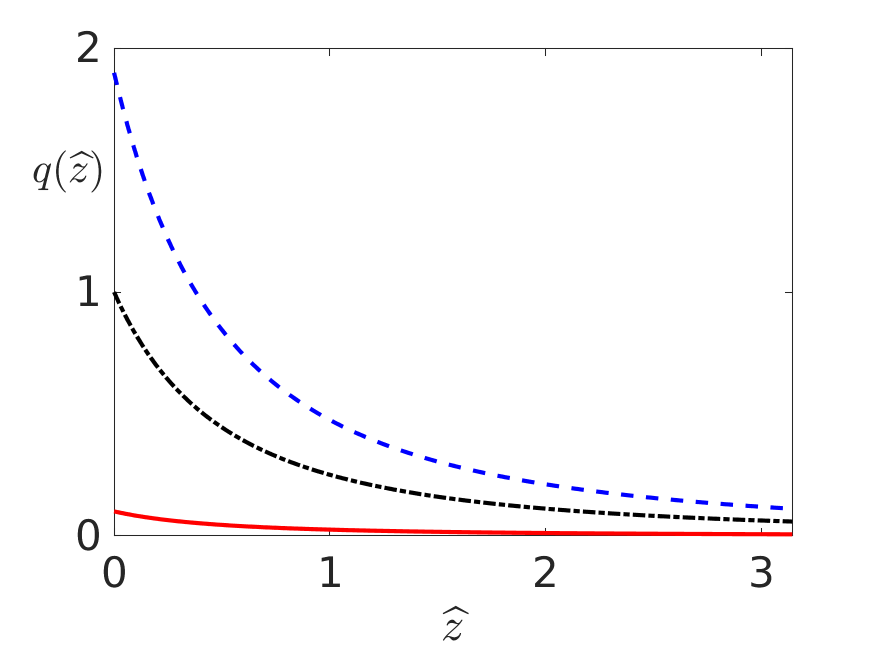}
\end{subfigure}
\caption{(\textbf{a}) Plots of the invariant function $q(\widehat{z})$ for $a = 1 = c$ in a logarithmic scale and several values of $b$: solid red ($b = 0.1$), dotted green ($b = 0.5$), dash-dotted black ($b = 1$), and dashed blue ($b = 2$). (\textbf{b}) Sketches of the invariant function $q(\widehat{z})$ for $a = 1 = b$ and several values of $c$: solid red ($c = 0.1$), dash-dotted black ($c = 1$), and dashed blue ($c = 1.9$).}	\label{invariant-plot}
\end{figure}

We start with the following definition of the generalized second PdHA problem.
\begin{Definition}
[Generalized second Paine--de Hoog--Anderssen (PdHA) problem]
The following Sturm--Liouville BVP in the Schr\"odinger (or Liouville normal) form is defined as a \emph{generalized second Paine--de Hoog--Anderssen (PdHA) problem} with the third-kind (Robin) boundary conditions:
\begin{align}
-\frac{d^2\widehat{y}}{d\widehat{z}^2} + q(\widehat{z}) \widehat{y} &= \lambda \widehat{y}, \qquad \qquad \widehat{z}_0 < \widehat{z} < \widehat{z}_1, \\
\alpha_0 \widehat{y}(\widehat{z}_0) - \alpha_1 \dot{\widehat{y}} (\widehat{z}_0) &= 0, \qquad \qquad
 \beta_0 \widehat{y}(\widehat{z}_1) +  \beta_1 \dot{\widehat{y}} (\widehat{z}_1) = 0,
\end{align}
where dot denotes the derivative with respect to $\widehat{z}$ and the invariant function $q$ is given by 
\begin{equation}
q(\widehat{z}) = \frac{c}{\left(a \widehat{z} + b \right)^n}, \qquad a, b, c > 0, \quad n \in \mathbb{N}.
\end{equation}
\end{Definition}
For $n = 2$, we call the above a \emph{modified second PdHA problem}. Figure~\ref{invariant-plot} displays plots of the invariant function $q(\widehat{z})$ for different values of parameters. The left panel of Figure~\ref{invariant-plot} sketches the function for fixed values of $a$ and $c$ but for different values of $b$ in a logarithmic scale. The right panel of Figure~\ref{invariant-plot} illustrates the plots of $q$ when $a$ and $b$ are fixed but the parameter $c$ is varied. In both instances, the function is decreasing for an increasing value \mbox{of $\widehat{z}$.}
\vspace{-6pt}
\end{paracol}
\nointerlineskip
\begin{paracol}{2}
\switchcolumn

For $n \neq 2$, the problem seems to be open. Our conjecture is that the spectrum exhibits a similar qualitative behaviour with the one for $n = 2$. Although for the BVP in the Liouville normal form can still be solved numerically, recovering the corresponding Sturm--Liouville problem in the canonical form turns out to be nontrivial.
\begin{Corollary}
In particular, for the modified second PdHA problem ($n = 2$), when $a = 1 = c$, $b = 0.1$, by imposing the Dirichlet conditions $\alpha_0 = 1 = \beta_0$, $\alpha_1 = 0 = \beta_1$, and, selecting the endpoints as $\widehat{z}_0 = 0$ and $\widehat{z}_1 = \pi$, the above BVP reduces to the well-known classical second PdHA problem~{\rm\cite{paine1981correction}}:
\begin{equation*}
-\frac{d^2\widehat{y}}{d\widehat{z}^2} + \frac{1}{\left(\widehat{z} + \frac{1}{10} \right)^2} \widehat{y} = \lambda \widehat{y}, \qquad \qquad 0 < \widehat{z} < \pi, \qquad \qquad \widehat{y}(0) = 0 = \widehat{y}(\pi).
\end{equation*}	
\end{Corollary}

We need the following theorem on a transformation of the Sturm--Liouville problem from the canonical form into the Schr\"odinger form~\cite{pryce1993numerical}:
\begin{Theorem}
The Sturm--Liouville problem in the canonical form with eigenvalue $\lambda$ and the corresponding eigenfunction $y(z)$ is given explicitly as follows:
\begin{align*}
-\frac{d}{dz} \left[u(z) \frac{dy}{dz} \right] = \lambda y,& \qquad \qquad z_0 < z < z_1, \\
a_0 y(z_0) - a_1 \frac{dy}{dz}(z_0) = 0,& \qquad \qquad b_0 y(z_1) + b_1 \frac{dy}{dz}(z_1) = 0.
\end{align*}
can be converted into the following BVP in the \emph{Schr\"odinger} (or \emph{Liouville normal}) \emph{form}
\begin{align*}
-\frac{d^2\widehat{y}}{d\widehat{z}^2} + q(\widehat{z}) \, \widehat{y} = \lambda \, \widehat{y},& \qquad \qquad \widehat{z}_0 < \widehat{z} < \widehat{z}_1,		\\
\widehat{\alpha}_0 \widehat{y}(\widehat{z_0}) - \widehat{\alpha}_1 \frac{d\widehat{y}}{d\widehat{z}}(\widehat{z_0}) = 0, & \qquad \qquad 
\widehat{ \beta}_0 \widehat{y}(\widehat{z_1}) + \widehat{ \beta}_1 \frac{d\widehat{y}}{d\widehat{z}}(\widehat{z_1}) = 0,	
\end{align*}
by performing \emph{Liouville's transformation}
\begin{equation}
\widehat{z} = \int \frac{dz}{\sqrt{u}}, \qquad v = u^{-1/4}, \qquad \text{and} \qquad y(z) = v(z) \, \widehat{y}(z),		\label{liou-trans}
\end{equation}
where
\begin{align*}
\alpha_0 &= a_0 v(\widehat{z}_0) - \frac{a_1}{\dot{z}(\widehat{z}_0)} \frac{dv}{d\widehat{z}}(\widehat{z}_0), \qquad \qquad 
&\alpha_1 = a_1 \frac{v(\widehat{z_0})}{\dot{z}(\widehat{z_0})}, \\
\beta_0 &= b_0 v(\widehat{z}_1) + \frac{b_1}{\dot{z}(\widehat{z_1})} \frac{dv}{d\widehat{z}}(\widehat{z}_1), \qquad \qquad 
&\beta_1 = b_1 \frac{v(\widehat{z_1})}{\dot{z}(\widehat{z_1})},
\end{align*}
and $q$ is the corresponding \emph{invariant function} of the ODE in the canonical form, given by
\begin{equation*}
q(\widehat{z}) = v \frac{d^2}{d\widehat{z}^2} \left(\frac{1}{v}\right).
\end{equation*}
\end{Theorem}

A detailed outline of the proof of this theorem is available in a recent preprint, where an application related to perturbed potential temperature distribution in an atmospheric boundary layer was also discussed via the WKB method and numerical simulation~\cite{karjanto2021perturbed}.

\begin{Proposition}
The following Sturm--Liouville problem in the canonical form with Dirichlet boundary conditions
\begin{align*}
-\frac{d}{dz}\left\{a\left[(1 - 2k)(z - d) \right]^{-4k/(1 - 2k)} \dfrac{dy}{dz} \right\} &= \lambda y, \qquad \qquad z_0 < z < z_1, \\
y(z_0) = 0 = y(z_1),
\end{align*}
where
\begin{align*}
   k &= \frac{1}{2}\left(-1 \pm \frac{1}{a}\sqrt{a^2 + 4c} \right), \qquad \qquad 
&z_0 &= \frac{b^{1 -2k}}{1 - 2k} + d, \\
   a &\neq 0, \qquad d \in \mathbb{R}, \qquad \qquad 
&z_1 &= \frac{(\pi + b)^{1 -2k}}{1 - 2k} + d,
\end{align*}
reduces to the modified second PdHA problem
\begin{equation}
-\frac{d^2\widehat{y}}{d\widehat{z}^2} + \frac{c}{(a\widehat{z} + b)^2} \widehat{y} = \lambda \widehat{y}, \qquad \qquad 0 < \widehat{z} < \pi, \qquad \qquad \widehat{y}(0) = 0 = \widehat{y}(\pi).
\end{equation}
\end{Proposition}
\begin{proof}
Since
\begin{equation*}
u(z) = \left[a(1 - 2k)(z - d) \right]^{-4k/(1 - 2k)},
\end{equation*}
then
\begin{equation*}
\widehat{z} = \int \frac{dz}{\sqrt{u}} = \frac{1}{a} \left\{ \left[a(1 - 2k)(z - d) \right]^{1/(1 - 2k)} - b \right\},
\end{equation*}
for which we can express $z$ in terms of $\widehat{z}$
\begin{equation*}
z = \frac{\left(a\widehat{z} + b\right)^{1 - 2k}}{a(1 - 2k)} + d.
\end{equation*}
We express the function $u$ in terms of the transformed variable $\widehat{z}$:
\begin{equation*}
u(\widehat{z}) = \left(a\widehat{z} + b\right)^{-4k},
\end{equation*} 
and the function $v$ using the transformation $v = u^{-1/4}$:
\begin{equation*}
v(\widehat{z}) = \left(a\widehat{z} + b\right)^{k}.
\end{equation*}
It follows that
\begin{align*}
\frac{d}{d\widehat{z}} \left(\frac{1}{v} \right) &= -a k \left(a\widehat{z} + b\right)^{-k - 1}, \\
\frac{d^2}{d\widehat{z}^2} \left(\frac{1}{v} \right) &= a^2 k(k + 1) \left(a\widehat{z} + b\right)^{-k - 2}, \\
q(\widehat{z}) = v \frac{d^2}{d\widehat{z}^2} \left(\frac{1}{v} \right) &= \frac{a^2 k(k + 1)}{\left(\widehat{z} + b\right)^2} = \frac{c}{\left(a\widehat{z} + b\right)^2}.
\end{align*}
The boundary points $\widehat{z}_0 = 0$ and $\widehat{z}_1 = \pi$ can be attained by substituting $z_0$ and $z_1$ into the expression for $\widehat{z}$, respectively. The proof is completed. \hfill {\qquad}
\end{proof}

\begin{Corollary}
The classical second PdHA problem can be achieved by taking $a = 1 = c$, $b = 0.1$, and $u(z) = \left[(2 \mp \sqrt{5})(z - d)\right]^{2(3 \pm \sqrt{5})}$.
\end{Corollary}

\begin{Example}
Consider a regular Sturm--Liouville problem with Dirichlet boundary conditions 
\begin{equation*}
-\frac{d}{dz} \left\{ \left[\left(2 + \sqrt{5} \right)(z - d)\right]^{2(3 - \sqrt{5})} \frac{dy}{dz} \right\} = \lambda y, \qquad \qquad y(0) = 0 = y(z_1). 
\end{equation*}
where $d = -(0.1)^{2 + \sqrt{5}}/(2 + \sqrt{5}) \approx 1.37 \times 10^{-5}$ and $z_1 = (\pi + 0.1)^{2 + \sqrt{5}}/(2 + \sqrt{5}) + d \approx 34.4068$.
The BVP in the Liouville normal form is given by
\begin{equation*}
-\frac{d^2\widehat{y}}{d\widehat{z}^2} + \frac{1}{\left(\widehat{z} + \frac{1}{10} \right)^2} \widehat{y} = \lambda \widehat{y}, \qquad 
\widehat{y}(0) = 0 = \widehat{y}(\pi).
\end{equation*}
\end{Example} 
Figure~\ref{compare-plot} shows the exponential behaviour in $z$ of the coefficient function $u(z)$ (left panel) and the corresponding invariant function $q(\widehat{z})$ in a logarithmic scale (right panel), which is considered in the classical second PdHA problem.

\clearpage
\end{paracol}
\nointerlineskip
\begin{figure}[H]
\widefigure
\begin{subfigure}{0.45\textwidth}
\caption{} \vspace{0.1cm}
\includegraphics[width=\linewidth]{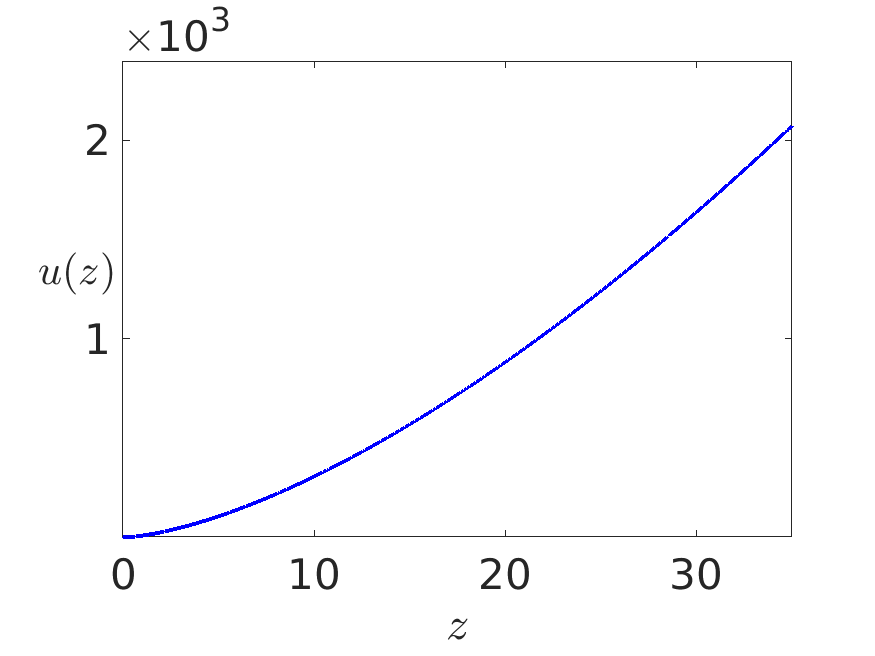}
\end{subfigure}	\hspace*{0.5cm}
\begin{subfigure}{0.45\textwidth}
\caption{} \vspace{0.1cm}
\includegraphics[width=\linewidth]{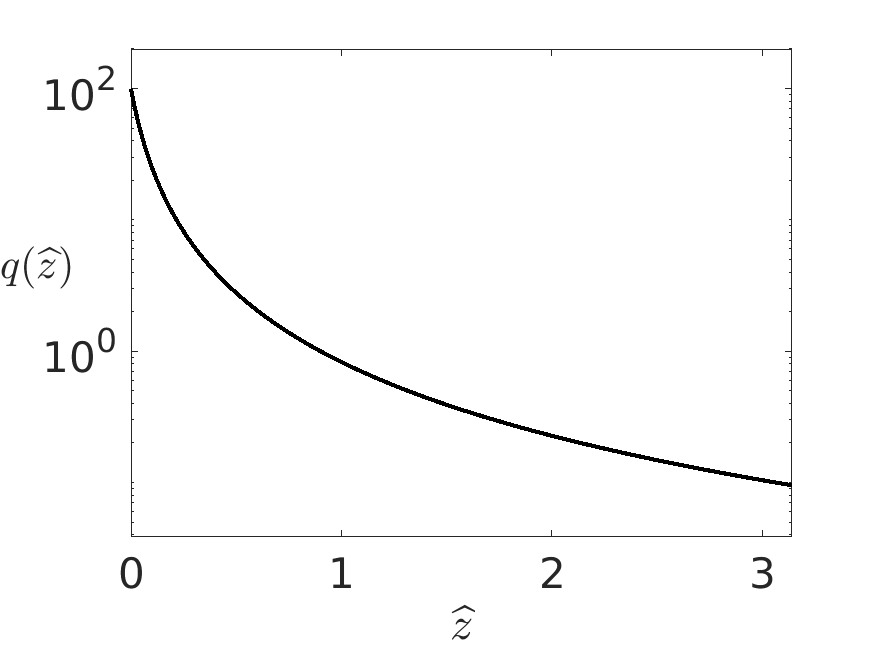}
\end{subfigure}
\caption{(\textbf{a}) A sketch of the coefficient function $u(z)$ showing an exponential behaviour as $z$ increases. (\textbf{b}) A plot of the invariant function $q(\widehat{z})$ in a logarithmic scale for the classical second PdHA problem. See the solid red curve in the left panel of Figure~\ref{invariant-plot}.}		\label{compare-plot}
\end{figure}
\begin{paracol}{2}
\switchcolumn

\section{Estimating the Smallest Eigenvalue}	\label{estimate}

Arnold et al. proposed a novel technique to approximate the smallest eigenvalue of a Sturm--Liouville problem by means of a localized landscape  and effective potential functions~\cite{arnold2019computing,filoche2016universal}. Let $\widehat{w}$ be a landscape function, then it satisfies the ODE
\begin{equation}
-\frac{d^2\widehat{w}}{d\widehat{z}^2} + q(\widehat{z}) \widehat{w} = 1, 
\end{equation}
with appropriate boundary conditions. An \emph{effective potential} $\widehat{W}$ is defined as $\widehat{W} = 1/\widehat{w}$, i.e., the reciprocal of the landscape function. Then, the lowest eigenvalue of the Sturm--Liouville problem is given by
\begin{equation}
\lambda_0 \approx \frac{5}{4} \widehat{W}_\text{min},
\end{equation}
where $\widehat{W}_\text{min}$ is the minimum value of the effective potential on the corresponding interval of the problem.

Without losing any generality, we can rewrite the invariant function $q$ in the following form:
\begin{equation*}
q(\widehat{z}) = \frac{\widehat{c}}{\left(\widehat{z} + \widehat{b} \right)^2}, \qquad \qquad \text{where} \qquad \widehat{c} = \frac{c}{a^2}, \qquad \text{and} \qquad \widehat{b} = \frac{b}{a}.
\end{equation*}
Hence, the corresponding landscape function for the modified second PdHA problem satisfies the following BVP:
\begin{equation*}
-\frac{d^2\widehat{w}}{d\widehat{z}^2} +\frac{\widehat{c}}{\left(\widehat{z} + \widehat{b} \right)^2} \widehat{w} = 1, \qquad \qquad  \widehat{w}(0) = 0 = \widehat{w}(\pi).
\end{equation*}
In what follows, we will pursue generalized conditions at the endpoints, i.e., the third-kind (Robin) boundary conditions, formulated as follows:
\begin{equation*}
\alpha_0 \widehat{w}(0) + \alpha_1 \frac{d\widehat{w}}{d\widehat{z}}(0) = w_0, \qquad \qquad
\beta_0 \widehat{w}(\pi) + \beta_1 \frac{d\widehat{w}}{d\widehat{z}}(\pi) = w_1.
\end{equation*} 
Observe that the first and second kinds or Dirichlet and Neumann boundary conditions can be acquired by taking $\alpha_1 = 0 = \beta_1$ and $\alpha_0 = 0 = \beta_0$, respectively.

Let $\widehat{w}(\widehat{z}) = (\widehat{z} + \widehat{b})^{\phi}$, $\phi \in \mathbb{R}$, be an Ansatz for the homogeneous ODE, and then, upon substitution, we obtain $\phi^2 - \phi - \widehat{c} = 0$, which solves
\begin{equation*}
\phi = \phi_{1,2} = \frac{1}{2} \left(1 \pm \sqrt{1 + 4 \widehat{c}} \right),
\end{equation*}
where the subscripts $1$ and $2$ correspond to the positive and negative signs, respectively. Observe that, for a special case of $\widehat{c} = 1$, the positive root $\phi_1 = \varphi$, the well-known golden ratio, and the conjugate root $\phi_2 = 1 - \varphi = -1/\varphi$, also known as the negative of the silver ratio. Thus, a complementary solution to the homogeneous ODE is given by
\begin{equation*}
\widehat{w}_c(\widehat{z}) = C_1 (\widehat{z} + \widehat{b})^{\phi_1} + C_2 (\widehat{z} + \widehat{b})^{\phi_2}, \qquad \qquad C_1, \; C_2 \in \mathbb{R}.
\end{equation*}
We seek the particular solution using the variation of parameters technique. Writing the Ansatz as 
\begin{equation*}
\widehat{w}_p(\widehat{z}) = r_1(\widehat{z}) (\widehat{z} + \widehat{b})^{\phi_1} + r_2(\widehat{z}) (\widehat{z} + \widehat{b})^{\phi_2},
\end{equation*}
we obtain
\begin{align*}
r_1(\widehat{z}) &= \int \frac{(\widehat{z} + \widehat{b})^{1 - \phi_1}}{(\phi_2 - \phi_1)} \, d\widehat{z} 
= \frac{(\widehat{z} + \widehat{b})^{2 - \phi_1}}{(\phi_2 - \phi_1)(2 - \phi_1)}, \\
r_2(\widehat{z}) &= \int \frac{-(\widehat{z} + \widehat{b})^{1 - \phi_2}}{(\phi_2 - \phi_1)} \, d\widehat{z} 
= \frac{-(\widehat{z} + \widehat{b})^{2 - \phi_2}}{(\phi_2 - \phi_1)(2 - \phi_2)}.
\end{align*}
Substituting these to the Ansatz for the particular solution yields
\begin{equation*}
\widehat{w}_p(\widehat{z}) = -\frac{(\widehat{z} + \widehat{b})^2}{(2 - \phi_1)(2 - \phi_2)} = \frac{(\widehat{z} + \widehat{b})^2}{\widehat{c} - 2}, \qquad \widehat{c} \neq 2.
\end{equation*}
Combining together both $\widehat{w}_c$ and $\widehat{w}_p$, we obtain a general solution for the landscape function:
\begin{equation*}
\widehat{w}(\widehat{z}) = C_1 (\widehat{z} + \widehat{b})^{\phi_1} + C_2 (\widehat{z} + \widehat{b})^{\phi_2} + \frac{(\widehat{z} + \widehat{b})^2}{\widehat{c} - 2}.
\end{equation*}
By imposing the Robin boundary conditions, the constant coefficients $C_1$ and $C_2$ satisfy the following matrix equation:
\end{paracol}
\begin{equation*}
\begin{bmatrix*}[c]
b^{\phi_1} \left(\alpha_0 + \frac{\alpha_1 \phi_1}{b} \right) &  b^{\phi_2} \left(\alpha_0 + \frac{\alpha_1 \phi_2}{b} \right) \\
(\pi + b)^{\phi_1} \left(\beta_0 + \frac{\beta_1 \phi_1}{\pi + b} \right) & (\pi + b)^{\phi_2} \left(\beta_0 + \frac{\beta_1 \phi_2}{\pi + b} \right)
\end{bmatrix*}
\begin{bmatrix*}[c]
C_1 \\ C_2
\end{bmatrix*} = 
\begin{bmatrix*}[c]
w_0 \\ w_1
\end{bmatrix*} + \frac{1}{2 - \widehat{c}}
\begin{bmatrix*}[c]
b(\alpha_0 b + 2\alpha_1) \\ (\pi + b) \left[\beta_0(\pi + b) + 2\beta_1 \right]
\end{bmatrix*},
\end{equation*}
\begin{paracol}{2}
\switchcolumn \noindent
which solves as $C_1 = C_1^R = \widehat{C}_1/\widehat{C}$ and $C_2 = C_2^R = \widehat{C}_2/\widehat{C}$, where
\end{paracol}
\begin{align*}
\widehat{C}_1 &= \left[w_0(2 - \widehat{c}) + b(\alpha_0 b + 2 \alpha_1) \right](\pi + b)^{\phi_2} \left(\beta_0 + \frac{\beta_1 \phi_2}{\pi + b} \right) 
- b^{\phi_2} \left(\alpha_0 + \frac{\alpha_1 \phi_2}{b} \right) \left\{w_1(2 - \widehat{c}) + (\pi + b)\left[\beta_0 (\pi + b) + 2 \beta_1 \right] \right\}, \\
\widehat{C}_2 &= b^{\phi_1} \left(\alpha_0 + \frac{\alpha_1 \phi_1}{b} \right) \left\{w_1(2 - \widehat{c}) + (\pi + b)\left[\beta_0 (\pi + b) + 2 \beta_1 \right] \right\}
- (\pi + b)^{\phi_1} \left(\beta_0 + \frac{\beta_1 \phi_1}{\pi + b} \right) \left[w_0(2 - \widehat{c}) + b(\alpha_0 b + 2 \alpha_1) \right], \\
\widehat{C} &= (2 - \widehat{c})\left\{b^{\phi_1} (\pi + b)^{\phi_2} \left(\alpha_0 + \frac{\alpha_1 \phi_1}{b} \right) \left(\beta_0 + \frac{\beta_1 \phi_2}{\pi + b} \right) - b^{\phi_2} (\pi + b)^{\phi_1} \left(\alpha_0 + \frac{\alpha_1 \phi_2}{b} \right)  \left(\beta_0 + \frac{\beta_1 \phi_1}{\pi + b} \right)\right\}.
\end{align*}
\begin{paracol}{2}
\switchcolumn 

The corresponding coefficients $C_1$ and $C_2$ for Dirichlet boundary conditions 
\begin{equation*}
\widehat{w}(0) = \widehat{w}_0 \qquad \text{and} \qquad
\widehat{w}(\pi) = \widehat{w}_1, \qquad \qquad \text{where} \qquad \widehat{w}_0 = \frac{w_0}{\alpha_0} \qquad \text{and} \qquad \widehat{w}_1 = \frac{w_1}{\beta_0},
\end{equation*}
are given as follows:
\begin{align*}
C_1 = C_1^{D} &= \frac{\left[\widehat{w}_0 (2 - \widehat{c}) + b^2 \right](\pi + b)^{\phi_2} - b^{\phi_2} \left[\widehat{w}_1(2 - \widehat{c}) + (\pi + b)^2\right]}{(2 - \widehat{c}) \left[ b^{\phi_1} (\pi + b)^{\phi_2} - b^{\phi_2} (\pi + b)^{\phi_1}\right]}, \\
C_2 = C_2^{D} &= \frac{b^{\phi_1} \left[\widehat{w}_1 (2 - \widehat{c}) + (\pi + b)^2 \right] - (\pi + b)^{\phi_1} \left[\widehat{w}_0(2 - \widehat{c}) + b^2\right]}{(2 - \widehat{c}) \left[ b^{\phi_1} (\pi + b)^{\phi_2} - b^{\phi_2} (\pi + b)^{\phi_1}\right]}.
\end{align*}
Similarly, we can also derive the coefficients for Neumann boundary conditions
\end{paracol}
\begin{equation*}
\frac{d\widehat{w}}{d\widehat{z}}(0) = \widetilde{w}_0 \qquad \text{and} \qquad
\frac{d\widehat{w}}{d\widehat{z}}(\pi) = \widetilde{w}_1, \qquad \qquad \text{where} \qquad \widetilde{w}_0 = \frac{w_0}{\alpha_1} \qquad \text{and} \qquad \widetilde{w}_1 = \frac{w_1}{\beta_1}.
\end{equation*}
\begin{paracol}{2}
\switchcolumn \noindent
The corresponding coefficients are given as follows:
\begin{align*}
C_1 = C_1^{N} &= \frac{\left[\widetilde{w}_0 (2 - \widehat{c}) + 2b \right](\pi + b)^{\phi_2 - 1} - b^{\phi_2 - 1} \left[\widetilde{w}_1(2 - \widehat{c}) + 2(\pi + b)\right]}{\phi_1 (2 - \widehat{c}) \left[ b^{\phi_1 - 1} (\pi + b)^{\phi_2 - 1} - b^{\phi_2 - 1} (\pi + b)^{\phi_1 - 1}\right]}, \\
C_2 = C_2^{N} &= \frac{b^{\phi_1 - 1} \left[\widetilde{w}_1 (2 - \widehat{c}) + 2(\pi + b) \right] - (\pi + b)^{\phi_1 - 1} \left[\widetilde{w}_0 (2 - \widehat{c}) + 2 b \right]}{\phi_2 (2 - \widehat{c}) \left[ b^{\phi_1 - 1} (\pi + b)^{\phi_2 - 1} - b^{\phi_2 - 1} (\pi + b)^{\phi_1 - 1} \right]}.
\end{align*}

The minimum value of effective potential, or the maximum value of the landscape function, occurs at $\widehat{z} = \widehat{z}_c$ such that $d\widehat{w}/d\widehat{z}(\widehat{z}_c) = 0$, i.e.,
\begin{equation*}
C_1 \phi_1 (\widehat{z}_c + \widehat{b})^{\phi_1 - 2} + C_2 \phi_2 (\widehat{z}_c + \widehat{b})^{\phi_2 - 2} = \frac{2}{(2 - \widehat{c})}.
\end{equation*}
We can solve this equation numerically using any standard technique of root-finding algorithms, e.g., the Newton--Raphson method or similar derivative-based techniques.

Figure~\ref{lfep-D} displays the plots of the landscape function $\widehat{w}(\widehat{z})$ and effective potential $\widehat{W}(\widehat{z})$ for several values of $\widehat{b}$ but for a chosen fixed value of $\widehat{c}$. The landscape function satisfies the BVP for the modified PdHA problem with homogeneous Dirichlet conditions at both endpoints. We observe that as $\widehat{b}$ increases, its slope at $\widehat{z} = 0$ also increases, which causes an increase in its maximum value. Consequently, the minimum values for the potential function decrease as $\widehat{b}$ increases, and an immediate consequence is also a decrease in the eigenvalue estimate. See the left panel of Figure~\ref{lambda-vs-b-c1}. 

\end{paracol}
\nointerlineskip
\begin{figure}[H]
\widefigure
\begin{subfigure}{0.45\textwidth}
\caption{} \vspace{0.1cm}
\includegraphics[width=\linewidth]{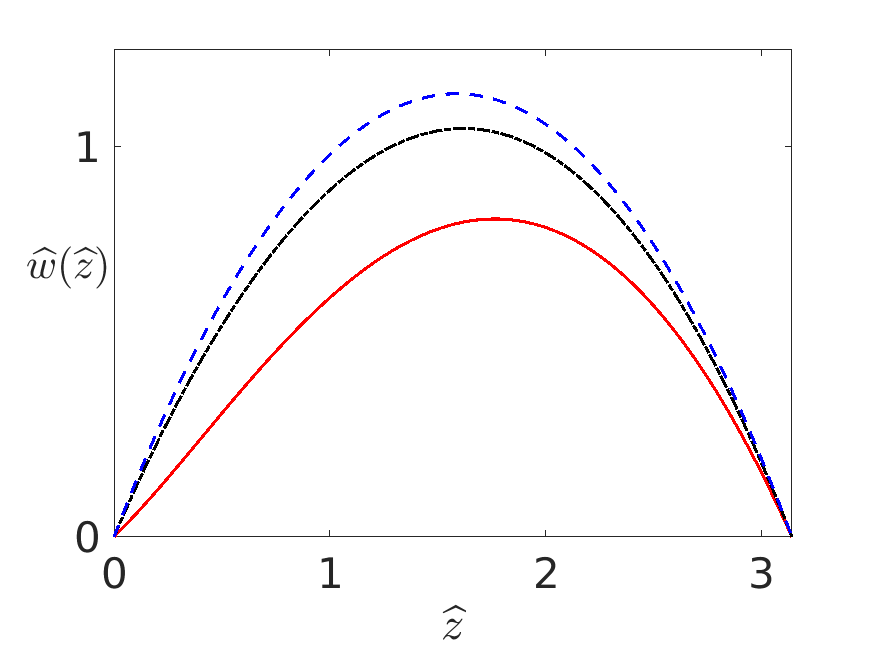}
\end{subfigure}	\hspace*{0.5cm}
\begin{subfigure}{0.45\textwidth}
\caption{} \vspace{0.1cm}
\includegraphics[width=\linewidth]{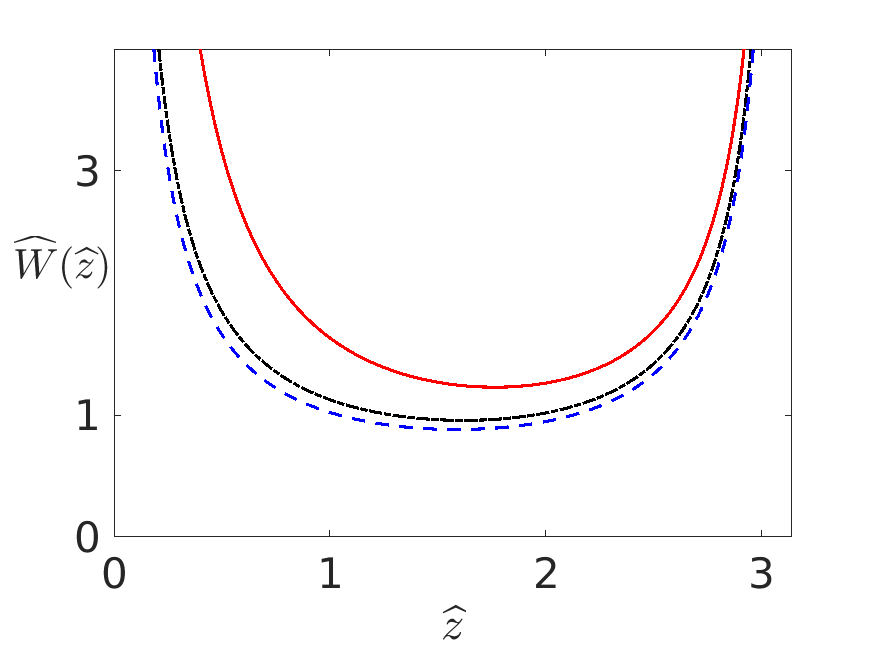}
\end{subfigure}
\caption{(\textbf{a}) A depiction of the landscape function $\widehat{w}(\widehat{z})$ that satisfies the Dirichlet boundary conditions $\widehat{w}(0) = 0 = \widehat{w}(\pi)$ for various values of $\widehat{b}$ and a fixed value of $\widehat{c} = 1$: solid red ($\widehat{b} = 0.1$), dash-dotted black ($\widehat{b} = 1$), and dashed blue ($\widehat{b} = 2$). (\textbf{b})~A sketch of the corresponding effective potential with the same parameter values as in the left panel.}	\label{lfep-D}
\end{figure}
\begin{paracol}{2}
\switchcolumn

Figure~\ref{lfep-N} displays the plots of the landscape function and the corresponding effective potential where the former satisfies the BVP for the modified second PdHA problem with nonhomogeneous Neumann boundary conditions. The slopes are chosen to be positive and negative unity at the left- and right-endpoints, respectively. Both panels depict the curves for a fixed value of $\widehat{c}$ and several values of $\widehat{b}$. Observe that although the slopes of the landscape function at both endpoints remain identical, the initial condition for each landscape function increases as the value of $\widehat{b}$ increases. Consequently, the maximum value also increases accordingly. Conversely, the effective potential exhibits a decreasing minimum value for an increasing value of $\widehat{b}$, which translates to a decreasing value of eigenvalue as $\widehat{b}$ is becoming large. See Figures~\ref{lambda-vs-b-c1} and~\ref{lambda-vs-b-c01-19}.

\end{paracol}
\nointerlineskip
\begin{figure}[H]
\widefigure
\begin{subfigure}{0.45\textwidth}
\caption{} \vspace{0.1cm}
\includegraphics[width=\linewidth]{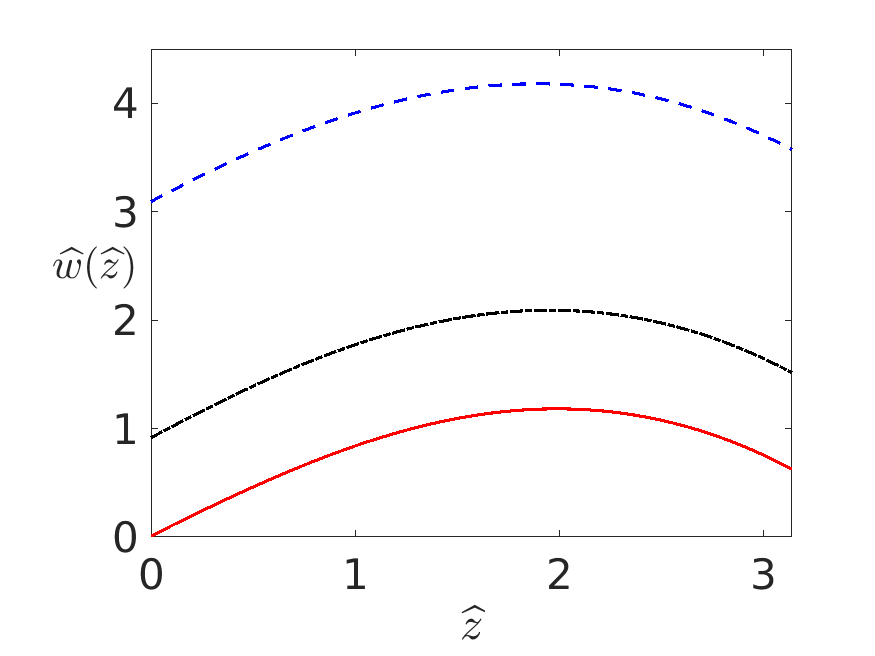}
\end{subfigure}	\hspace*{0.5cm}
\begin{subfigure}{0.45\textwidth}
\caption{} \vspace{0.1cm}
\includegraphics[width=\linewidth]{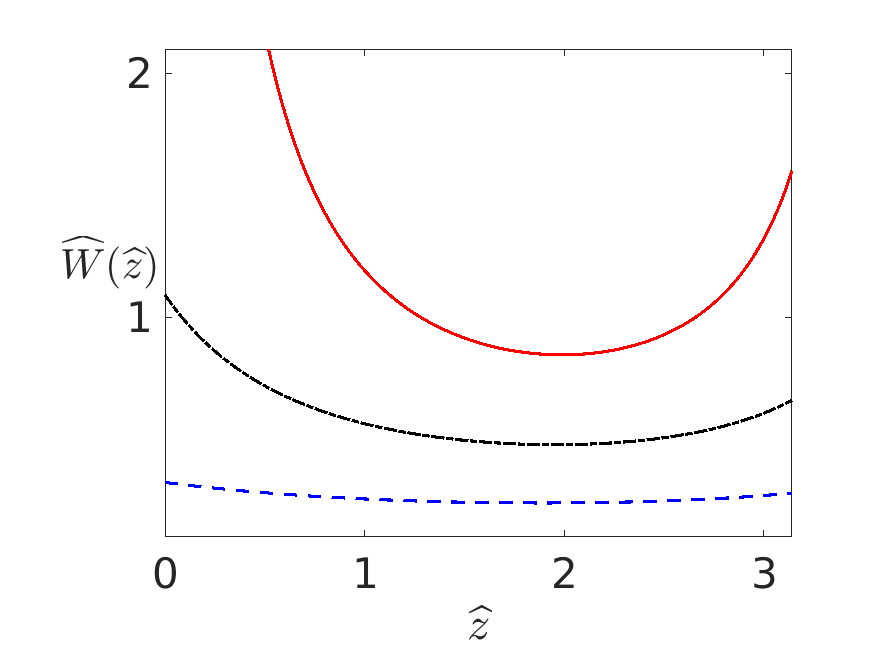}
\end{subfigure}
\caption{Similar to Figure~\ref{lfep-D} but for Neumann boundary conditions $\widetilde{w}_0 = 1 = -\widetilde{w}_1$. The solid red curves now correspond to $\widehat{b} = 0.3$.}	\label{lfep-N}
\end{figure}
\begin{paracol}{2}
\switchcolumn

\vspace{-12pt}

\section{Numerical Comparison}		\label{compare}

In this section, we compare the lowest-order eigenvalues obtained through the landscape function with the ones from numerical simulation. We utilize \textsc{MATLAB}'s finite difference code \emph{bvp4c} for the latter, accompanied by call functions of the system of the ODE, boundary conditions, and initial guess to solve the corresponding BVP with an unknown parameter.  

\subsection{Eigenvalue Comparison}

In this subsection, we provide comparisons of the eigenvalues of the modified second PdHA problem between the ones obtained from the landscape function and numerical simulations. We also examine further where the prescribed boundary conditions are either Dirichlet or Neumann type. 

Figure~\ref{lambda-vs-b-c1} displays the plots of the lowest-order eigenvalue $\lambda_0$ as a function of the parameter $\widehat{b}$. The left panel of Figure~\ref{lambda-vs-b-c1} depicts a comparison between the eigenvalues obtained from the landscape function and the results from numerical simulations for several values of parameter $\widehat{c}$ and particular Dirichlet boundary conditions. Although the approximation overestimates the numerical calculations, the eigenvalue generally exhibits remarkable quantitative agreement for a wide range of parameter values. 

The right panel of Figure~\ref{lambda-vs-b-c1} presents a comparison of the eigenvalues obtained by an estimate from the landscape function and numerical simulations for particular Neumann boundary conditions with $\widehat{c} = 1$. We observe that although the qualitative behavior for both outcomes is similar, i.e., generally decreasing for increasing value of $\widehat{b}$, its quantitative behavior is not. For $\widehat{b} \leq 2$, the approximation is far from accurate, whilst for $b > 2$, the eigenvalue estimate demonstrates a better approximation with its numerical counterpart despite the former still overestimating the latter. 

Figure~\ref{lambda-vs-b-c01-19} shows a similar comparison with the right panel of Figure~\ref{lambda-vs-b-c1} but for a small value of $\widehat{c}$ (left panel) and for the value of $\widehat{c}$ where the landscape function is nearly singular (right panel). A distinct qualitative behavior emerges from this finding. For small $\widehat{c}$, the eigenvalue estimate has a better quantitative agreement with the numerical results. As $\widehat{c}$ increases, it turns out that the accuracy fades away for smaller values of $\widehat{b}$, as we also observed in the right panel of Figure~\ref{lambda-vs-b-c1}. For this particular case, we also noticed that for $\widehat{b} > 2$, the estimate improves significantly whereas it was completely inaccurate for $\widehat{b} \leq 2$.

\end{paracol}
\nointerlineskip
\begin{figure}[H]
\widefigure
\begin{subfigure}{0.45\textwidth}
\caption{} \vspace{0.1cm}
\includegraphics[width=\linewidth]{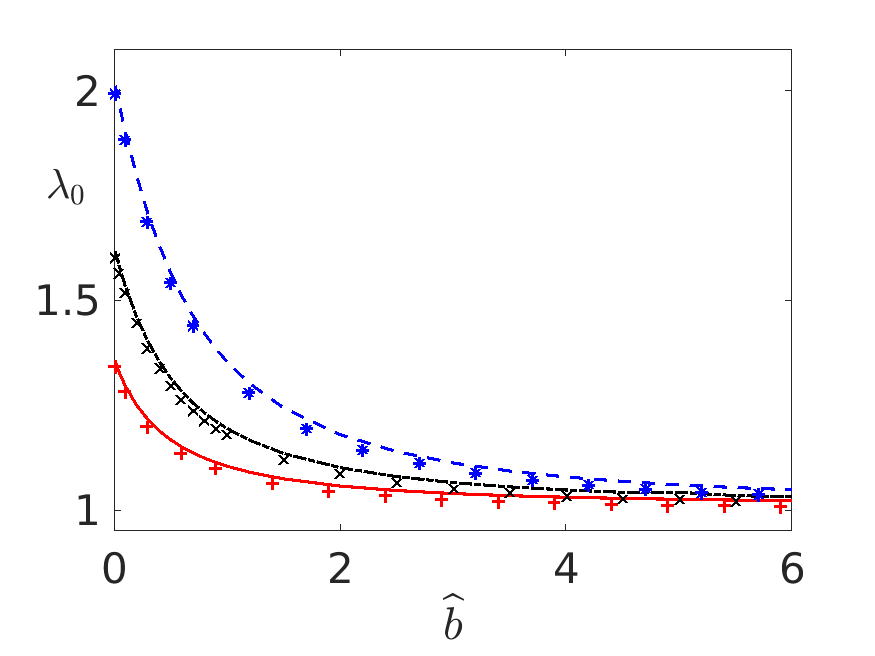}
\end{subfigure}	\hspace*{0.5cm}
\begin{subfigure}{0.45\textwidth}
\caption{} \vspace{0.1cm}
\includegraphics[width=\linewidth]{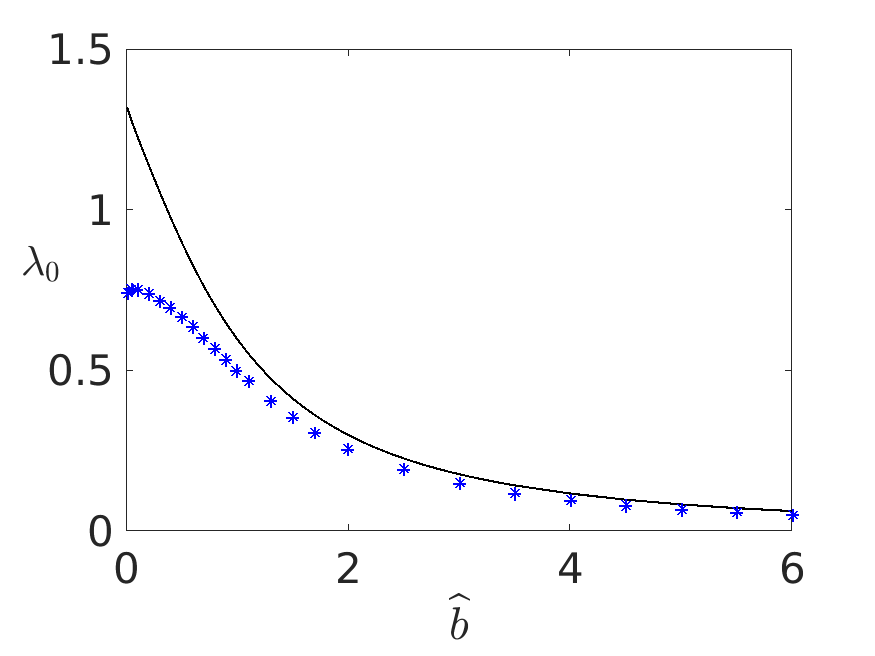}
\end{subfigure}
\caption{(\textbf{a}) Plot of the lowest-order eigenvalue $\lambda_0$ as a function of $\widehat{b}$ for several values of $\widehat{c}$ for Dirichlet boundary conditions $\widehat{w}_0 = 0 = \widehat{w}_1$: solid red ($\widehat{c} = 0.5$), dashed-dotted black ($\widehat{c} = 1$), and dashed blue ($\widehat{c} = 1.9$). The corresponding red pluses, black crosses, and blue asterisks are the eigenvalues obtained numerically. (\textbf{b}) A depiction of the lowest-order eigenvalue $\lambda_0$ as a function of $\widehat{b}$ for Neumann boundary conditions $\widetilde{w}_0 = 1 = -\widetilde{w}_1$ where $\widehat{c} = 1$. The black solid curve is an eigenvalue estimate obtained from the landscape function, whereas the blue asterisks correspond to the eigenvalues obtained by numerical simulations.}		\label{lambda-vs-b-c1}
\end{figure}

\begin{figure}[H]
\widefigure
\begin{subfigure}{0.45\textwidth}
\caption{} \vspace{0.1cm}
\includegraphics[width=\linewidth]{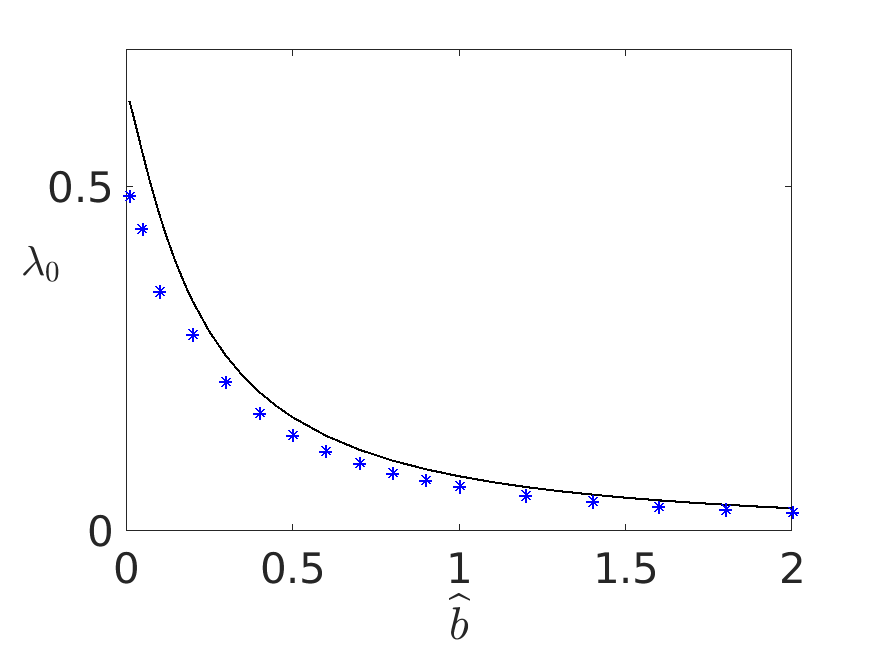}
\end{subfigure}	\hspace*{0.5cm}
\begin{subfigure}{0.45\textwidth}
\caption{} \vspace{0.1cm}
\includegraphics[width=\linewidth]{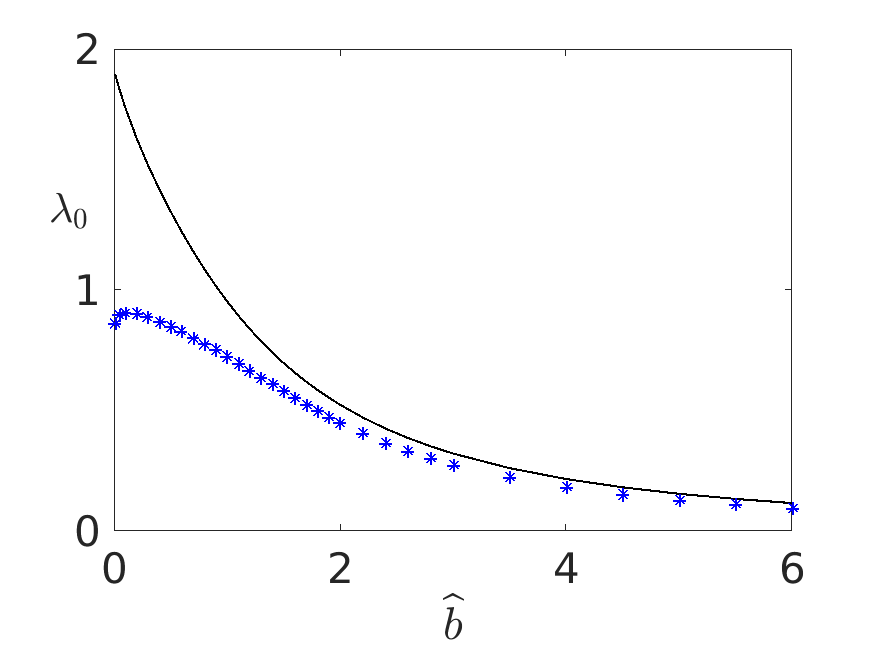}
\end{subfigure}
\caption{Similar to the right panel of Figure~\ref{lambda-vs-b-c1} but for different values of $\widehat{c}$: (\textbf{a}) $\widehat{c} = 0.1$, (\textbf{b}) $\widehat{c} = 1.9$. In all cases, the eigenvalues obtained from the landscape function and effective potential overestimates the eigenvalues obtained numerically. For $\widehat{c} = 1.9$, the estimate is worse for $\widehat{b} \leq 2$ but it becomes better as the value of $\widehat{b}$ increases, particularly for $\widehat{b} > 2$.}		\label{lambda-vs-b-c01-19}
\end{figure}
\begin{figure}[H]
\widefigure
\begin{subfigure}{0.45\textwidth}
\caption{} \vspace{0.1cm}
\includegraphics[width=\linewidth]{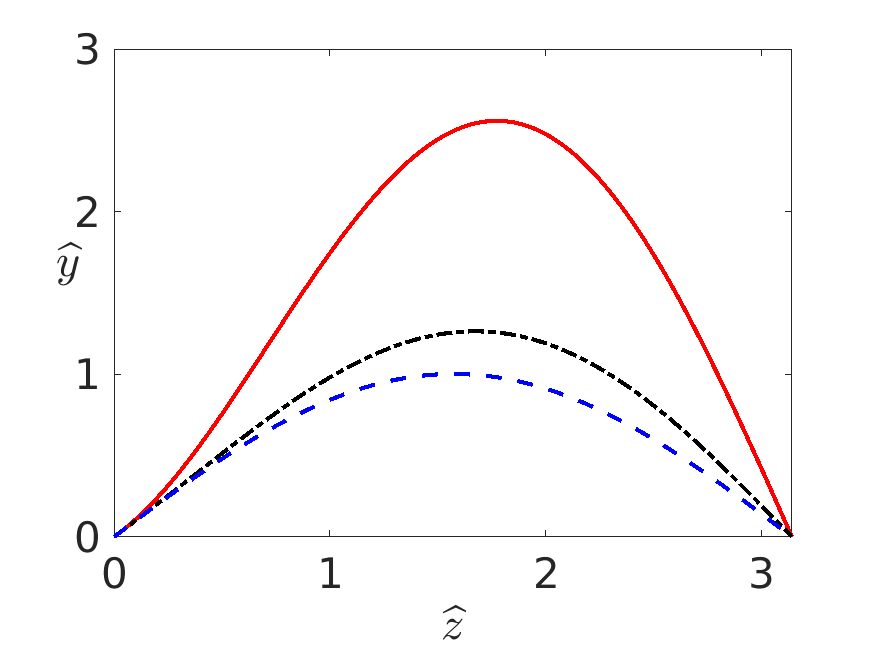}
\end{subfigure}	\hspace*{0.5cm}
\begin{subfigure}{0.45\textwidth}
\caption{} \vspace{0.1cm}
\includegraphics[width=\linewidth]{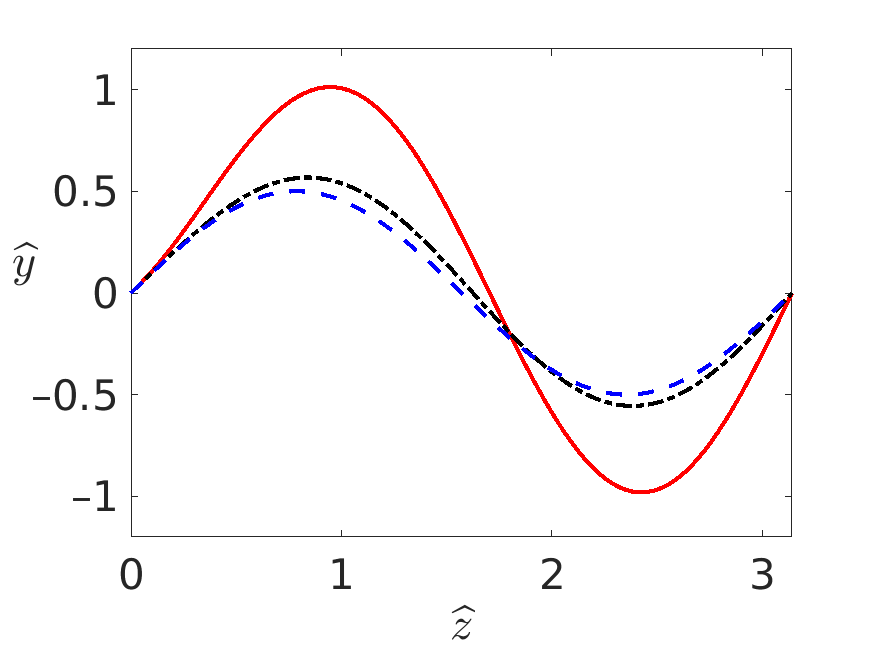}
\end{subfigure}
\caption{Plots of eigenfunctions satisfying Dirichlet boundary conditions $\widehat{y}(0) = 0 = \widehat{y}(\pi)$ obtained numerically for a fixed value of $\widehat{c} = 1$ and several values of $\widehat{b}$: solid red ($\widehat{b} = 0.1$), dash-dotted black ($\widehat{b} = 0.5$), and dashed blue ($\widehat{b} = 6$). (\textbf{a}) Plots of the lowest-order eigenfunctions with eigenvalues $\lambda_0 = 1.520$, $1.297$, and $1.018$, respectively. (\textbf{b}) Plots of the second-order eigenfunctions with eigenvalues $\lambda_1 = 4.493$, $4.416$, and $4.018$, respectively.}		\label{efD-numer}
\end{figure}
\begin{figure}[H]
\widefigure
\begin{subfigure}{0.45\textwidth}
\caption{} \vspace{0.1cm}
\includegraphics[width=\linewidth]{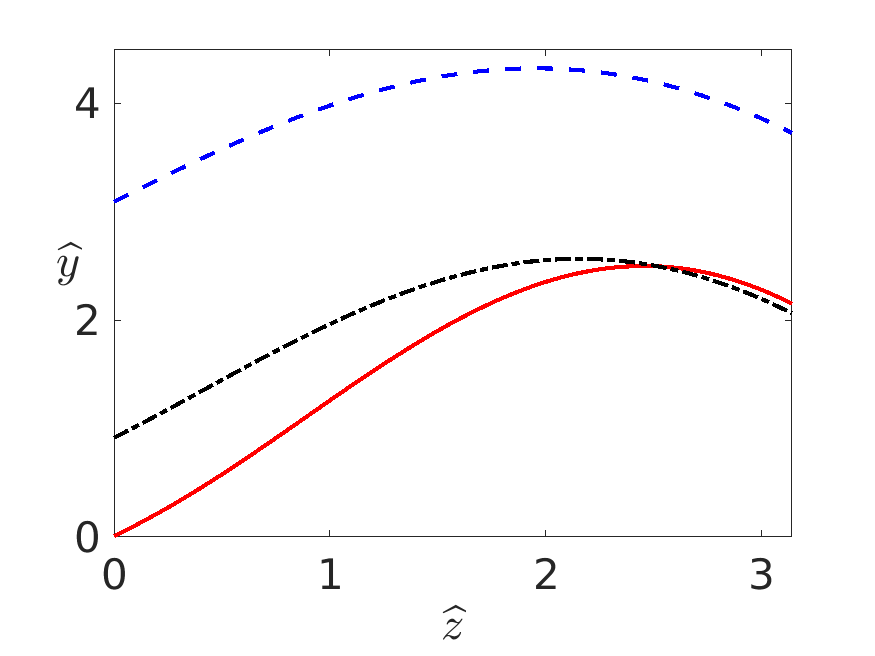}
\end{subfigure}	\hspace*{0.5cm}
\begin{subfigure}{0.45\textwidth}
\caption{} \vspace{0.1cm}
\includegraphics[width=\linewidth]{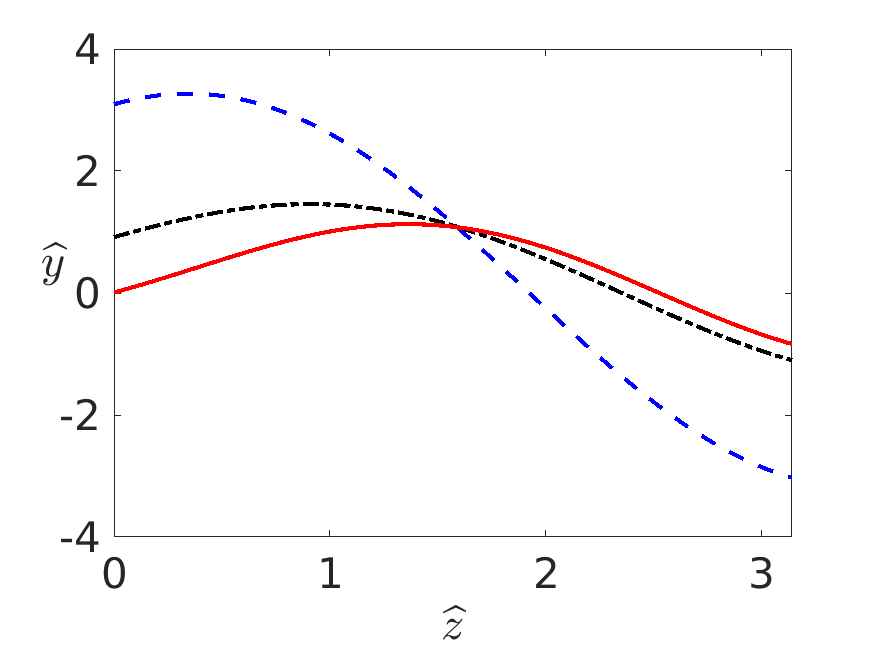}
\end{subfigure}
\caption{Similar as in Figure~\ref{efD-numer} but now with the eigenfunctions satisfy Neumann boundary conditions $\dot{\widehat{y}}(0) = 1 = -\dot{\widehat{y}}(\pi)$. The solid red curves correspond to $\widehat{b} = 0.3$. (\textbf{a}) Plots of the lowest-order eigenfunctions with eigenvalues $\lambda_0 = 0.718$, $0.500$, and $0.254$, respectively. (\textbf{b}) Plots of the second-order eigenfunctions with eigenvalues $\lambda_1 = 2.082$, $1.371$, and $1.119$, respectively.}		\label{efN-numer}
\end{figure}

\begin{paracol}{2}
\switchcolumn

\subsection{Numerical Eigenfunction}

In this subsection, we present the corresponding eigenfunctions of the modified second PdHA BVP for Dirichlet and Neumann conditions. We only focus on the first two lowest-order eigenfunctions since higher-order eigenfunctions can be calculated accordingly by simply modifying the values of the initial guess for the eigenvalue and/or eigenfunction.

Figure~\ref{efD-numer} displays the eigenfunctions for the modified second PdHA problem with particular Dirichlet boundary conditions. The lowest-order and the second-order eigenfunctions are presented in the left and right panels of Figure~\ref{efD-numer}, respectively. A fixed value of $\widehat{c}$ is prescribed and different curves correspond to various values of $\widehat{b}$. We observe that for increasing { values} of $\widehat{b}$, both maximum and minimum values of the eigenfunctions decrease. A similar trend also occurs for the eigenvalues, as we  discussed in the previous subsection. From these two lowest-order eigenfunctions, even though an increase in the value of $\widehat{b}$ is rather large, e.g., for $\widehat{b} \geq 1$, remarkably, the eigenfunction profiles do not change drastically. The opposite situation occurs when the value of $\widehat{b}$ is relatively small, i.e., $\widehat{b} \leq 0.5$, the decrease in the maximum value and quantitative profile of the eigenfunctions is quite dramatic. 
\end{paracol}
\nointerlineskip
\begin{paracol}{2}
\switchcolumn

Figure~\ref{efN-numer} illustrates similar eigenfunctions for the modified second PdHA BVP but with specified Neumann conditions. The left and right panels correspond to the first- and second-order eigenfunction profiles for different values of $\widehat{b}$ and a fixed parameter $\widehat{c}$. Compared to the previous case, we recognize that the eigenfunctions exhibit completely different quantitative and qualitative behaviours when we modify the boundary conditions from Dirichlet to Neumann type. While the slopes of the eigenfunction remain identical at both endpoints, the initial values do not. Rather, they rise in tandem with an increase of $\widehat{b}$. The boundary values at the right endpoints, however, do not necessarily follow this trend. An increase in the maximum values of the eigenfunctions is consistent with the decrease in eigenvalues as we  confirmed earlier. Different from the previous case, a relatively large change in the parameter value, e.g., from $\widehat{b} = 0.5$ to $6.0$, appears to correspond with a substantial change in the eigenfunction profiles quantitatively.

\clearpage
\end{paracol}
\nointerlineskip
\begin{paracol}{2}
\switchcolumn

\section{Conclusions}	\label{conclude}

We  considered a special case of a regular Sturm--Liouville BVP in this article. By transforming the problem into the Liouville normal form, we focus on the second PdHA problem. Although the classical PdHA problem {was} introduced four decades ago, a modified version of the second PdHA problem appears to be absent from the literature. The exposition above suggests that the modified second PdHA problem might provide additional insights when it comes to the behavior of the eigenvalue and corresponding eigenfunctions with respect to changes in the parameter values.

By incorporating the landscape function and effective potential, we can estimate the lowest-order eigenvalue without solving the corresponding eigenvalue problem.  We find that the approach is pretty remarkable even though numerical techniques could deliver the desired output within seconds. In particular, for special cases of the corresponding invariant function that corresponds to the modified second PdHA problem, both landscape and effective potential functions can be solved and expressed analytically. 

By prescribing Dirichlet- and Neumann-type boundary conditions, we  demonstrated the comparison of eigenvalues between the estimate and numerical simulations. Generally, the eigenvalues obtained from the landscape function overestimate the numerical results quantitatively but they exhibit relatively excellent qualitative agreement. The only exception occurred when we imposed Neumann boundary conditions and selected small values of parameter $\widehat{b}$ but for larger values of $\widehat{c}$.

\vspace{16pt}

\authorcontributions{Conceptualization, Methodology, Software, Validation, Formal Analysis, Investigation, Resources, Data Curation, Writing -- Original Draft Preparation, Writing -- Review and Editing, Visualization, Project Administration, N.K.}

\funding{This research received no external funding.}

\institutionalreview{Not applicable.}

\informedconsent{Not applicable.}

\dataavailability{The data presented in this study are available on request from the corresponding author.} 

\acknowledgments{The author gratefully acknowledges all anonymous reviewers for their useful and constructive feedback.}

\conflictsofinterest{The author declares no conflict of interest.}

\end{paracol}

\section*{Dedication}
\noindent
The author would like to dedicate this article to his late father Zakaria Karjanto (Khouw~Kim~Soey, 許金瑞) who introduced and taught him the alphabet, numbers and the calendar in his early childhood. Karjanto senior was born in Tasikmalaya, West Java, Japanese-occupied Dutch~East~Indies on 1~January~1944 (Saturday~Pahing) and died in Bandung, West Java, Indonesia on 18~April~2021 (Sunday~Wage). 

\reftitle{References}

\end{document}